\documentclass[11pt]{article}

\usepackage[T1]{fontenc}
\usepackage[utf8]{inputenc}
\usepackage[english]{babel}

\usepackage{amsmath}
\usepackage{amssymb}
\usepackage{amsthm}
\usepackage{mathtools}

\usepackage[a4paper,margin=30mm]{geometry}
\usepackage{microtype}
\usepackage[hidelinks]{hyperref}

\theoremstyle{plain}
\newtheorem{theorem}{Theorem}[section]
\newtheorem{proposition}[theorem]{Proposition}
\newtheorem{corollary}[theorem]{Corollary}
\newtheorem{lemma}[theorem]{Lemma}

\theoremstyle{definition}
\newtheorem{definition}[theorem]{Definition}

\theoremstyle{remark}
\newtheorem{remark}[theorem]{Remark}

\newcommand{\K}{\mathbb{K}}

\newcommand{\Q}{\mathbb{Q}}

\newcommand{\PP}{\mathbb{P}}

\newcommand{\cO}{\mathcal{O}}
\newcommand{\cI}{\mathcal{I}}
\newcommand{\cH}{\mathcal{H}}
\newcommand{\cE}{\mathcal{E}}
\newcommand{\cT}{\mathcal{T}}
\newcommand{\cQ}{\mathcal{Q}}

\DeclareMathOperator{\Supp}{Supp}
\DeclareMathOperator{\Bl}{Bl}
\DeclareMathOperator{\rk}{rk}
\DeclareMathOperator{\length}{length}
\DeclareMathOperator{\mult}{mult}

\DeclareMathOperator{\vdim}{vdim}
\DeclareMathOperator{\edim}{edim}

\title{Fat base schemes and unexpected hypersurfaces of type $(d+k,d)$}

\author{Grzegorz Malara}

\begin{document}

\maketitle

\begin{abstract}
We study hypersurfaces of degree $d+k$ containing a prescribed zero-dimensional fat point scheme $Z\subset\PP^n$ and having multiplicity at least $d$ along a general codimension-two linear subspace. We associate with this interpolation problem a vector bundle on $\PP^1$ whose splitting type determines the dimensions of all linear systems of type $(d+k,d)$.

For a fat point of multiplicity $m$, we determine the exact number of conditions detected by this bundle and introduce the corresponding effective length. This leads to an excess formula separating the contributions of the splitting type, the part of the fat structure invisible to degree-$k$ fibrewise sections, and the Hilbert-function defect of the fixed scheme.

We also prove a truncation theorem showing that multiplicities larger than $k$ contribute only a uniform twist of the bundle and, geometrically, fixed hyperplane components. Finally, we exhibit explicit planar and three-dimensional examples showing that genuinely unexpected phenomena persist for nonreduced fixed base schemes, including geometrically irreducible unexpected hypersurfaces.
\end{abstract}

\section{Introduction}

Unexpected hypersurfaces arise from interpolation problems in which a prescribed fixed scheme and a general multiple linear subspace fail to impose the expected number of independent conditions. In one of the standard forms of the problem, one fixes a finite set of points $Z\subset\PP^n$ and asks for hypersurfaces of degree $d+k$ containing $Z$ and vanishing to order at least $d$ along a general codimension-two linear subspace. Thus one studies the spaces \[ [I_Z\cap I_{\cH}^d]_{d+k}, \] and compares their actual dimensions with the dimensions predicted by naive interpolation.

The systematic study of unexpected curves originates in the work of Cook, Harbourne, Migliore and Nagel \cite{CHMN}, motivated in part by earlier results of Di Gennaro, Ilardi and Vall\`es \cite{DGV} and by the relation between interpolation and splitting types of logarithmic bundles observed by Faenzi and Vall\`es \cite{FV}. The theory was subsequently developed in several directions, including higher dimensional unexpected hypersurfaces \cite{BMSS,HMNT,HMTG}, special line arrangements \cite{DMO,MPTG}, determinantal and interpolation-matrix methods \cite{Matrixwise,MatrixwiseReloaded}, and structural questions concerning the persistence of unexpectedness \cite{FGHM}; see also the survey \cite{HMNSurvey}. For systems of type $(d+k,d)$, a syzygy-bundle approach was developed in the plane in \cite{MTG} and extended to higher-dimensional projective spaces in \cite{JMTG}. Most of the general theory in these works is formulated for reduced fixed configurations, although nonreduced fixed schemes occur in particular examples.

Nonreduced schemes have appeared in a different role in the theory. In particular, Szpond \cite{SzpondFat} constructed unexpected hypersurfaces with several general fat points, using Fermat-type configurations as the fixed reduced base locus. Unexpected curves also remain closely connected with combinatorial questions for line arrangements. Recent work of Dimca and Pokora \cite{DP} relates the generic splitting type, unexpected curves, and a generalized form of Terao's conjecture.

In the present paper we consider the complementary situation: the moving codimension-two centre remains general, whereas the prescribed fixed scheme \( Z=m_1P_1+\cdots+m_sP_s \) is itself allowed to be nonreduced. Our purpose is to develop the splitting-bundle formalism in this setting and to determine precisely which part of the fat structure is detected by it.

Let $\cH\subset\PP^n$ be a general codimension-two linear subspace. Blowing up $\PP^n$ along $\cH$ resolves the projection from $\cH$ and gives a morphism \( \rho:\Bl_{\cH}\PP^n\longrightarrow\PP^1. \) For every $k\geq0$ we associate with $Z$ a vector bundle \( \cE_{k,Z,\cH} = \rho_*\bigl(\cI_Z\otimes\cO_Y(kL)\bigr) \) on $\PP^1$. The basic observation is that, for every $d\geq0$, 
\[ [I_Z\cap I_{\cH}^d]_{d+k} \simeq H^0\bigl(\PP^1,\cE_{k,Z,\cH}(d)\bigr). \]
 Consequently, if 
\[ \cE_{k,Z,\cH} \simeq \bigoplus_{i=1}^{R_k}\cO_{\PP^1}(-a_i), \]
 then the splitting type $(a_1,\ldots,a_{R_k})$ simultaneously determines the dimensions of all systems of type $(d+k,d)$ associated with the fixed scheme $Z$.

The first new feature of the nonreduced setting is that the bundle does not necessarily detect the whole scheme-theoretic length of a fat point. For a point of multiplicity $m$ we determine the precise number
\[
\lambda_{n,k}(m) = \sum_{b=0}^{\min\{k,m-1\}} (m-b)\binom{n+b-2}{n-2}
\]
of conditions visible to degree-$k$ fibrewise sections. This leads to the effective length
\[
\Lambda_k(Z) = \sum_i\lambda_{n,k}(m_i)
\]
and to the jet defect
\(\delta_k(Z) = \deg Z-\Lambda_k(Z).\)
In particular, the sum of the splitting exponents is governed by the effective length rather than by the ordinary length:
\[
\sum_{i=1}^{R_k}a_i = \Lambda_k(Z)-B_k.
\]
This distinction leads to an excess formula which separates three different sources of speciality. If
\[
T = \sum_{i=1}^{R_k}(a_i-d-1)_+
\]
is the splitting tail and \( \eta_Z(d+k) = h^1\bigl(\PP^n,\cI_Z(d+k)\bigr) \) is the Hilbert-function defect of the fixed scheme, then
\[
\dim[I_Z\cap I_{\cH}^d]_{d+k} - \vdim_{\mathrm{base}} = T+\delta_k(Z)-\eta_Z(d+k).
\]
Thus the discrepancy between actual and virtual dimensions decomposes into a contribution coming from the splitting type, a contribution coming from the part of the fat structure invisible to the bundle, and a correction measuring the speciality of $Z$ itself.

A second feature of the theory is that sufficiently large fixed multiplicities have a simple geometric interpretation. We prove that multiplicities larger than $k$ may be truncated to $k$ at the cost of a uniform twist of the bundle. More precisely, if \( Z_{\leq k} = \sum_i\min\{m_i,k\}P_i \) and \( r=\sum_i(m_i-k)_+, \) then
\[
\cE_{k,Z,\cH} \simeq \cE_{k,Z_{\leq k},\cH}(-r).
\]
Geometrically, the excess multiplicities give fixed hyperplane components through $\cH$. Consequently, genuinely new phenomena for fat base schemes may be sought among schemes whose fixed multiplicities are at most $k$. In particular, for $k=1$ the problem reduces completely to the reduced support, whereas for $k=2$ double points provide the first essentially nonreduced case.

The final section illustrates the theory on explicit configurations. For a configuration of thirteen points in $\PP^2$ we identify a distinguished triple \( P_1, P_3, P_{11} \) with two related unexpected behaviours. Doubling these three points and taking $k=2$ gives splitting type \( (4,6,6) \) and produces a unique unexpected sextic with a general quadruple point. For a general point $Q$ this sextic is unique up to scalar and geometrically irreducible. Tripling the same three points and taking $k=3$ gives splitting type \( (5,5,5,7) \) and a three-dimensional system of octics with a general quintuple point, while the expected dimension is two. This system again contains absolutely irreducible members. The relation between the two examples is made explicit: a two-dimensional subsystem of the octics is obtained by multiplying the unexpected sextic by conics through the three distinguished points and the general point, while the excess produces an additional independent direction.

We also reconsider a three-dimensional Fermat configuration in which one coordinate point is replaced by a triple point. For $k=3$ the associated splitting type is \[ (1,2,3,3,3,3,3,4,4,4). \] It gives both an unexpected quartic containing a general line and an unexpected three-dimensional system of quintics vanishing doubly along a general line. The latter contains geometrically irreducible members, showing that the nonreduced fixed scheme produces phenomena which are not explained merely by fixed components or reducible constructions.

The paper is organized as follows. In Section~\ref{sec:fat-bundle} we construct the vector bundle associated with a fat point scheme and identify its sections with the relevant interpolation spaces. Section~\ref{sec:ambient-effective} determines the ambient bundle and introduces the effective length. In Section~\ref{sec:splitting-excess} we derive the splitting formula and the excess formula. Section~\ref{sec:reduction} is devoted to the truncation theorem and its interpretation in terms of fixed hyperplane components. Finally, Section~\ref{sec:examples} contains the planar and three-dimensional examples described above.

\section{The bundle associated with a fat point scheme}
\label{sec:fat-bundle}

Let $\K$ be an algebraically closed field of characteristic zero and let \( S=\K[x_0,\ldots,x_n] \) be the homogeneous coordinate ring of $\PP^n$.

Throughout this section, let \( Z=m_1P_1+\cdots+m_sP_s\subset\PP^n \) be a zero-dimensional fat point scheme, with saturated homogeneous ideal
\[ I_Z=\bigcap_{i=1}^s I(P_i)^{m_i}. \] 
Let $\cH\subset\PP^n$ be a codimension-two linear subspace satisfying
\begin{enumerate}
\item[(H1)] $\cH\cap\Supp(Z)=\emptyset$,
\item[(H2)] the points $q_i=\rho(P_i)\in\PP^1$ are pairwise distinct, i.e.\ $\cH$ meets no line joining two points of $\Supp(Z)$.
\end{enumerate}
All results of Sections \ref{sec:fat-bundle}, \ref{sec:ambient-effective} and \ref{sec:reduction} hold for every $\cH$ satisfying (H1) and (H2). In Section \ref{sec:splitting-excess} we moreover assume
\begin{enumerate}
\item[(H3)] $\cH$ belongs to the nonempty open subset of the Grassmannian on which the splitting type of $E_{k,Z,\cH}$ is constant (see Section \ref{subsec:splitting-type}).
\end{enumerate}

For integers $d,k\geq0$, we are interested in the vector space \( \left[I_Z\cap I_{\cH}^d\right]_{d+k}. \) It consists of hypersurfaces of degree $d+k$ containing the fat scheme $Z$ and vanishing along $\cH$ with multiplicity at least $d$. Our aim is to interpret this space as a space of sections of a vector bundle on $\PP^1$.

\subsection{Projection from a codimension-two linear space}

Let \( \pi:Y=\Bl_{\cH}\PP^n\longrightarrow\PP^n \) be the blow-up of $\PP^n$ along $\cH$. We denote by $E$ the exceptional divisor and put \( L=\pi^*\cO_{\PP^n}(1). \) After a change of coordinates we may assume that \( \cH=V(x_0,x_1). \) Projection from $\cH$ is the rational map
\[
\PP^n\dashrightarrow\PP^1, \qquad [x_0:\cdots:x_n]\longmapsto[x_0:x_1].
\]
The blow-up resolves this rational map and gives a morphism \( \rho:Y\longrightarrow\PP^1. \) The pencil of hyperplanes through $\cH$ has divisor class $L-E$ on $Y$. Consequently,
\begin{equation}
\label{eq:rho-pullback}
\rho^*\cO_{\PP^1}(1) \simeq \cO_Y(L-E).
\end{equation}
The fibre of $\rho$ over a point of $\PP^1$ is the strict transform of a hyperplane containing $\cH$. Since $\cH$ is a Cartier divisor inside such a hyperplane, blowing up the hyperplane along $\cH$ does not change it. Hence every fibre of $\rho$ is isomorphic to $\PP^{n-1}$.

Moreover, \( \cO_Y(L)|_{\rho^{-1}(q)} \simeq \cO_{\PP^{n-1}}(1) \) for every $q\in\PP^1$.

\begin{definition}
\label{def:Tk}
For $k\geq0$ we define \( \cT_k:=\rho_*\cO_Y(kL). \)
\end{definition}

Since the fibres of $\rho$ are isomorphic to $\PP^{n-1}$ and \( H^i\bigl(\PP^{n-1},\cO_{\PP^{n-1}}(k)\bigr)=0 \), for \(i>0,\) and \( k\geq0,\) cohomology and base change show that $\cT_k$ is a vector bundle of rank
\begin{equation*}
R_k := \rk\cT_k = h^0\bigl(\PP^{n-1},\cO_{\PP^{n-1}}(k)\bigr) = \binom{n+k-1}{n-1}.
\end{equation*}

The splitting of $\cT_k$ will be determined explicitly in the next section.

\subsection{The bundle defined by the fat point scheme}

Since the support of $Z$ is disjoint from $\cH$, the blow-up morphism $\pi$ is an isomorphism in a neighbourhood of $Z$. We therefore use the same symbol $Z$ for its inverse image in $Y$, and denote its ideal sheaf on $Y$ by $\cI_Z$.

\begin{definition}
\label{def:EkZ}
For $k\geq0$ we define \( \cE_{k,Z,\cH} := \rho_*\bigl(\cI_Z\otimes\cO_Y(kL)\bigr). \)
\end{definition}

The inclusion \( \cI_Z\otimes\cO_Y(kL) \hookrightarrow \cO_Y(kL) \) induces an injection \( \cE_{k,Z,\cH} \hookrightarrow \cT_k. \) Indeed, tensoring the defining exact sequence of $Z$ by $\cO_Y(kL)$ gives
\[
0 \longrightarrow \cI_Z\otimes\cO_Y(kL) \longrightarrow \cO_Y(kL) \longrightarrow \cO_Z(kL) \longrightarrow 0.
\]
Applying $\rho_*$ yields
\[
0 \longrightarrow \cE_{k,Z,\cH} \longrightarrow \cT_k \longrightarrow \rho_*\cO_Z(kL).
\]
In particular, the quotient \( \cT_k/\cE_{k,Z,\cH} \) is supported on finitely many points of $\PP^1$. Since $\cE_{k,Z,\cH}$ is a torsion-free coherent sheaf on the smooth curve $\PP^1$, it is locally free. Moreover,
\begin{equation*}
\rk\cE_{k,Z,\cH} = \rk\cT_k = R_k.
\end{equation*}
Thus the fat point scheme $Z$ determines, for every $k\geq0$, a vector bundle of rank $R_k$ on $\PP^1$.

The relevance of this bundle to unexpected hypersurfaces is contained in the following observation.

\begin{proposition}
\label{prop:sections-identification}
For all integers $d,k\geq0$ there is a natural isomorphism \begin{equation*} \left[I_Z\cap I_{\cH}^d\right]_{d+k} \simeq H^0\bigl( \PP^1, \cE_{k,Z,\cH}(d) \bigr).
\end{equation*}
\end{proposition}

\begin{proof}
A form of degree $d+k$ vanishing along $\cH$ with multiplicity at least $d$ corresponds, after pulling back to the blow-up $Y$, to a section of \( \cO_Y\bigl((d+k)L-dE\bigr). \) Since $Z$ is disjoint from $\cH$, imposing in addition vanishing along $Z$ gives
\[
\left[I_Z\cap I_{\cH}^d\right]_{d+k} \simeq H^0\left( Y, \cI_Z\otimes \cO_Y\bigl((d+k)L-dE\bigr) \right).
\]
Indeed, away from $Z$ the additional ideal $\cI_Z$ is trivial, while near $Z$ the morphism $\pi$ is an isomorphism and the exceptional divisor is absent. Thus the two vanishing conditions are imposed independently on the two disjoint loci.

Using \eqref{eq:rho-pullback}, we have \( (d+k)L-dE = kL+d(L-E), \) and hence
\[
\cO_Y\bigl((d+k)L-dE\bigr) \simeq \cO_Y(kL)\otimes\rho^*\cO_{\PP^1}(d).
\]
Therefore
\[
\begin{aligned}
\left[I_Z\cap I_{\cH}^d\right]_{d+k}
&\simeq H^0\left( Y, \cI_Z\otimes\cO_Y(kL) \otimes\rho^*\cO_{\PP^1}(d) \right) \\
 &= H^0\left( \PP^1, \rho_*\left( \cI_Z\otimes\cO_Y(kL) \otimes\rho^*\cO_{\PP^1}(d) \right) \right).
\end{aligned}
\]
By the projection formula,
\[
\rho_*\left( \cI_Z\otimes\cO_Y(kL) \otimes\rho^*\cO_{\PP^1}(d) \right) \simeq \cE_{k,Z,\cH}\otimes\cO_{\PP^1}(d),
\]
which proves the assertion.
\end{proof}

\begin{remark}
\label{rem:no-R1}
No vanishing of a higher direct image is needed in Proposition~\ref{prop:sections-identification}. Indeed, for every coherent sheaf $\mathcal F$ on $Y$ one has \( H^0(Y,\mathcal F) = H^0(\PP^1,\rho_*\mathcal F), \) and the projection formula applies directly to the push-forward.
\end{remark}

By the Birkhoff--Grothendieck theorem there exist integers \( a_1,\ldots,a_{R_k} \) such that \( \cE_{k,Z,\cH} \simeq \bigoplus_{i=1}^{R_k} \cO_{\PP^1}(-a_i). \) At this stage, however, neither the sum of the integers $a_i$ nor their relation with the multiplicities of the points of $Z$ is apparent. The next section is devoted to the structure of the ambient bundle $\cT_k$ and to the amount of information about a fat point that is detected by $\cT_k$.

\section{The ambient bundle and the effective length}
\label{sec:ambient-effective}

We now determine the structure of the ambient bundle $\cT_k$ and measure the contribution of a fat point to the quotient $\cT_k/\cE_{k,Z,\cH}$. The latter need not coincide with the full length of the fat point: for fixed $k$, sufficiently high jets are not detected by sections of $\cO_Y(kL)$.

\subsection{The ambient bundle}

The morphism \( \rho:Y\longrightarrow\PP^1 \) admits a standard description as a projective bundle. We use the convention that $\PP(\mathcal V)$ parametrizes one-dimensional subspaces of the fibres of $\mathcal V$.

\begin{proposition}
\label{prop:Y-projective-bundle}
There is an isomorphism
\[
Y \simeq \PP_{\PP^1} \bigl( \cO_{\PP^1}(-1)\oplus \cO_{\PP^1}^{\oplus(n-1)} \bigr)
\]
over $\PP^1$, under which \( \cO_Y(L) \simeq \cO_{\PP(\mathcal V)}(1), \) where \( \mathcal V = \cO_{\PP^1}(-1)\oplus \cO_{\PP^1}^{\oplus(n-1)}. \)
\end{proposition}

\begin{proof}
Recall that, after a change of coordinates, \( \cH=V(x_0,x_1) \) and the morphism $\rho$ resolves the projection \[ [x_0:\cdots:x_n] \longmapsto [x_0:x_1]. \] Let $[s:t]\in\PP^1$. The fibre of \( \mathcal V = \cO_{\PP^1}(-1)\oplus \cO_{\PP^1}^{\oplus(n-1)} \) over $[s:t]$ is naturally identified with \( \langle(s,t)\rangle \oplus \K^{n-1}. \)
Thus a point of $\PP(\mathcal V)$ over $[s:t]$ can be represented as
\[
[\lambda(s,t):v_2:\cdots:v_n].
\]
The corresponding point of $\PP^n$ is \( [\lambda s:\lambda t:v_2:\cdots:v_n]. \) This realizes $\PP(\mathcal V)$ as the resolution of the projection from $\cH$, hence as the blow-up $Y=\Bl_{\cH}\PP^n$. Under this identification, the relative hyperplane bundle $\cO_{\PP(\mathcal V)}(1)$ is the pull-back of $\cO_{\PP^n}(1)$. Hence it coincides with $\cO_Y(L)$.
\end{proof}

The projective bundle description gives the splitting of $\cT_k$
immediately.

\begin{proposition}
\label{prop:Tk-splitting}
For every $k\geq0$ one has
\begin{equation}
\label{eq:Tk-splitting}
	\cT_k \simeq \bigoplus_{j=0}^k \cO_{\PP^1}(j)^{ \oplus\binom{n+k-j-2}{n-2} }.
\end{equation}
In particular,
\begin{equation*}
\rk\cT_k = R_k = \binom{n+k-1}{n-1}, \qquad \deg\cT_k = B_k := \binom{n+k-1}{n}.
\end{equation*}
\end{proposition}

\begin{proof}
With the notation of Proposition~\ref{prop:Y-projective-bundle},
\[
\cT_k = \rho_*\cO_Y(kL) \simeq \rho_*\cO_{\PP(\mathcal V)}(k) \simeq \operatorname{Sym}^k(\mathcal V^\vee).
\]
Since \( \mathcal V^\vee \simeq \cO_{\PP^1}(1) \oplus \cO_{\PP^1}^{\oplus(n-1)}, \) we obtain
\[
\operatorname{Sym}^k(\mathcal V^\vee) \simeq \bigoplus_{j=0}^k \cO_{\PP^1}(j) \otimes \operatorname{Sym}^{k-j} \bigl( \cO_{\PP^1}^{\oplus(n-1)} \bigr).
\]
Now
\[
\rk \operatorname{Sym}^{k-j} \bigl( \cO_{\PP^1}^{\oplus(n-1)} \bigr) = \binom{n+k-j-2}{n-2},
\]
which gives \eqref{eq:Tk-splitting}.

Summing the multiplicities gives
\[
\rk\cT_k = \sum_{j=0}^k \binom{n+k-j-2}{n-2} = \binom{n+k-1}{n-1}.
\]
Similarly,
\[
\deg\cT_k = \sum_{j=0}^k j\binom{n+k-j-2}{n-2} = \binom{n+k-1}{n},
\]
where the last equality is the standard binomial identity.
\end{proof}

\subsection{The effective contribution of a fat point}

Let \( P\in\PP^n\setminus\cH \) and put \( q=\rho(P)\in\PP^1. \) We first determine the contribution of the fat point $mP$ to the quotient of $\cT_k$.

\begin{definition}
\label{def:lambda-nkm}
For integers $n\geq2$, $k\geq0$, and $m\geq1$, define
\begin{equation}
\label{eq:lambda-definition}
\lambda_{n,k}(m) := \sum_{b=0}^{\min\{k,m-1\}} (m-b) \binom{n+b-2}{n-2}.
\end{equation}
We call $\lambda_{n,k}(m)$ the \emph{$k$-effective length} of the fat point $mP$.
\end{definition}

The terminology is justified by the following local computation.

\begin{proposition}
\label{prop:local-effective-length}
Let
\[
Q_{k,mP} := \operatorname{im} \left( \cT_k \longrightarrow \rho_*\cO_{mP}(kL) \right).
\]
Then $Q_{k,mP}$ is a torsion sheaf supported at $q=\rho(P)$ and \( \length Q_{k,mP} = \lambda_{n,k}(m). \)
\end{proposition}

\begin{proof}
Since $\rho$ is a $\PP^{n-1}$-bundle, we may choose a local parameter $t$ on $\PP^1$ at $q$ and local fibre coordinates \( u_1,\ldots,u_{n-1} \) centred at $P$. Thus, after completion, \( \widehat{\cO}_{Y,P} \simeq \K[[t,u_1,\ldots,u_{n-1}]] \) and
\[
\cO_{mP} \simeq \K[[t,u_1,\ldots,u_{n-1}]] / (t,u_1,\ldots,u_{n-1})^m.
\]
On a neighbourhood of $q$, the bundle $\cT_k$ has a basis given, after dehomogenizing on the fibre, by the monomials
\[
u^\beta = u_1^{\beta_1}\cdots u_{n-1}^{\beta_{n-1}}, \qquad |\beta|\leq k.
\]
Hence the image in $\cO_{mP}$ is spanned by the monomials \( t^a u^\beta \) satisfying \( |\beta|\leq k \) and \( a+|\beta|\leq m-1. \) 
Fix $b=|\beta|$. There are \( \binom{n+b-2}{n-2} \) monomials $u^\beta$ of total degree $b$ in $n-1$ variables, while for each such monomial the exponent $a$ may take the values \( 0,1,\ldots,m-1-b. \) Thus there are $m-b$ possible values of $a$. Summing over \( 0\leq b\leq\min\{k,m-1\} \) gives
\[
\length Q_{k,mP} = \sum_{b=0}^{\min\{k,m-1\}} (m-b) \binom{n+b-2}{n-2},
\]
as claimed.
\end{proof}

The ordinary length of the fat point $mP\subset\PP^n$ is
\[
\length(mP) = \binom{n+m-1}{n}.
\]
The effective length agrees with the ordinary length precisely as long as the fibre degree $k$ is large enough to detect all the monomials occurring in the local algebra.

\begin{corollary}
\label{cor:lambda-full-length}
If \( k\geq m-1, \) then
\begin{equation*}
\lambda_{n,k}(m) = \binom{n+m-1}{n} = \length(mP).
\end{equation*}
\end{corollary}

\begin{proof}
If $k\geq m-1$, the condition $|\beta|\leq k$ is automatic for every monomial surviving modulo \( (t,u_1,\ldots,u_{n-1})^m. \) Thus the monomials counted in Proposition~\ref{prop:local-effective-length} form a basis of the entire local algebra of $mP$.
\end{proof}

In particular, a discrepancy between the effective length and the ordinary length can occur only when \( m\geq k+2. \)

We now return to the fat point scheme \( Z=m_1P_1+\cdots+m_sP_s. \) Since $\cH$ is chosen general, we may and will assume in addition that the points \[ q_i:=\rho(P_i)\in\PP^1 \] are pairwise distinct. Equivalently, $\cH$ avoids the finitely many secant lines joining two distinct points of $\Supp(Z)$.

\begin{definition}
\label{def:Lambda-delta}
The \emph{$k$-effective length} of $Z$ is
\begin{equation}
\label{eq:Lambda-definition}
\Lambda_k(Z) := \sum_{i=1}^s \lambda_{n,k}(m_i).
\end{equation}
We also define the \emph{$k$-jet defect} of $Z$ by
\begin{equation*}
\delta_k(Z) := \deg Z-\Lambda_k(Z), \end{equation*} 
where \( \deg Z = \sum_{i=1}^s \binom{n+m_i-1}{n}. \)
\end{definition}

By construction, \( \delta_k(Z)\geq0. \) Moreover, \( \delta_k(Z)=0 \) whenever \( m_i\leq k+1\), for every \(i\).

Consider now the morphism \( \cT_k \longrightarrow \rho_*\cO_Z(kL) \) and denote its image by \( \cQ_{k,Z,\cH}. \) Then the definition of $\cE_{k,Z,\cH}$ gives an exact sequence
\begin{equation}
\label{eq:E-T-Q}
0 \longrightarrow \cE_{k,Z,\cH} \longrightarrow \cT_k \longrightarrow \cQ_{k,Z,\cH} \longrightarrow 0.
\end{equation}

\begin{proposition}
\label{prop:Q-total-length}
The torsion sheaf $\cQ_{k,Z,\cH}$ has length
\begin{equation*}
\length\cQ_{k,Z,\cH} = \Lambda_k(Z).
\end{equation*}
Consequently, \( \deg\cE_{k,Z,\cH} = B_k-\Lambda_k(Z). \)
\end{proposition}

\begin{proof}
Since the points $q_i=\rho(P_i)$ are pairwise distinct, the local contributions of the fat points $m_iP_i$ are supported at distinct points of $\PP^1$. Therefore
\[
\length\cQ_{k,Z,\cH} = \sum_{i=1}^s \length Q_{k,m_iP_i}.
\]
By Proposition~\ref{prop:local-effective-length}, this equals \( \sum_{i=1}^s \lambda_{n,k}(m_i) = \Lambda_k(Z). \)

Finally, since \eqref{eq:E-T-Q} is an exact sequence whose first two terms are vector bundles of the same rank and whose quotient is torsion, one has
\[
\deg\cT_k = \deg\cE_{k,Z,\cH} + \length\cQ_{k,Z,\cH}.
\]
Using Proposition~\ref{prop:Tk-splitting} gives \( \deg\cE_{k,Z,\cH} = B_k-\Lambda_k(Z). \)
\end{proof}

\begin{remark}
\label{rem:not-full-rhoOZ}
In general, the morphism \( \cT_k \longrightarrow \rho_*\cO_Z(kL) \) need not be surjective. Its image has length $\Lambda_k(Z)$, whereas \( \length\cO_Z = \deg Z. \) Thus the difference \( \delta_k(Z) = \deg Z-\Lambda_k(Z) \) measures exactly the part of the fat point structure which is invisible to the degree-$k$ fibrewise sections encoded by $\cT_k$.
\end{remark}

\section{Splitting type and the excess formula}
\label{sec:splitting-excess}

We now combine the construction of Section~\ref{sec:fat-bundle} with the effective-length computation of Section~\ref{sec:ambient-effective}. This gives a numerical description of the linear systems \( \left[I_Z\cap I_{\cH}^d\right]_{d+k} \) in terms of the splitting type of the bundle $\cE_{k,Z,\cH}$.

\subsection{Splitting type}\label{subsec:splitting-type}

By the Birkhoff--Grothendieck theorem, there exist integers \( a_1\leq\cdots\leq a_{R_k} \) such that
\begin{equation}
\label{eq:E-splitting}
\cE_{k,Z,\cH} \simeq \bigoplus_{i=1}^{R_k} \cO_{\PP^1}(-a_i).
\end{equation}

\begin{definition}
\label{def:splitting-type}
The ordered tuple \( (a_1,\ldots,a_{R_k}) \) will be called the \emph{$k$-th splitting type of $Z$ with respect to $\cH$}.
\end{definition}

For a general choice of $\cH$, this splitting type is constant on a nonempty open subset of the corresponding Grassmannian. In what follows, $\cH$ is always assumed to belong to this open subset.

The first numerical constraint on the splitting type follows immediately from the computation of the degree of $\cE_{k,Z,\cH}$.

\begin{proposition}
\label{prop:sum-exponents}
The splitting exponents satisfy
\begin{equation}
\label{eq:sum-exponents}
\sum_{i=1}^{R_k}a_i = \Lambda_k(Z)-B_k,
\end{equation}
where \( B_k=\binom{n+k-1}{n}. \)
\end{proposition}

\begin{proof}
By Proposition~\ref{prop:Q-total-length}, \( \deg\cE_{k,Z,\cH} = B_k-\Lambda_k(Z). \) On the other hand, the splitting \eqref{eq:E-splitting} gives
\[
\deg\cE_{k,Z,\cH} = -\sum_{i=1}^{R_k}a_i.
\]
Comparing the two expressions proves the assertion.
\end{proof}

Combining the splitting with Proposition~\ref{prop:sections-identification} gives the basic dimension formula.

\begin{theorem}
\label{thm:dimension-splitting}
Let \( \cE_{k,Z,\cH} \simeq \bigoplus_{i=1}^{R_k} \cO_{\PP^1}(-a_i). \) Then, for every $d\geq0$,
\begin{equation*}
\dim \left[I_Z\cap I_{\cH}^d\right]_{d+k} = \sum_{i=1}^{R_k} \max\{0,d-a_i+1\}.
\end{equation*}
\end{theorem}

\begin{proof}
By Proposition~\ref{prop:sections-identification}, \( \left[I_Z\cap I_{\cH}^d\right]_{d+k} \simeq H^0\bigl( \PP^1,\cE_{k,Z,\cH}(d) \bigr). \) Using \eqref{eq:E-splitting}, we get \( \cE_{k,Z,\cH}(d) \simeq \bigoplus_{i=1}^{R_k} \cO_{\PP^1}(d-a_i). \) 
Since \( h^0\bigl(\PP^1,\cO_{\PP^1}(r)\bigr) = \max\{0,r+1\}, \) the result follows.
\end{proof}

\subsection{Virtual and expected dimensions}

We now compare the dimension in Theorem~\ref{thm:dimension-splitting} with the expected dimension of the corresponding interpolation problem.

Put \( D=d+k \) and let
\[
N_D := h^0\bigl(\PP^n,\cO_{\PP^n}(D)\bigr) = \binom{D+n}{n}.
\]

We first record the number of conditions imposed on degree-$D$ hypersurfaces by requiring multiplicity at least $d$ along a codimension-two linear subspace.

\begin{definition}
\label{def:Cndk}
We set
\begin{equation*}
C_{n,d,k} := \dim \left( S_D/ [I_{\cH}^d]_D \right).
\end{equation*}
Thus $C_{n,d,k}$ is the codimension of $[I_{\cH}^d]_D$ in $S_D$.
\end{definition}

After choosing coordinates such that \( I_{\cH}=(x_0,x_1), \) a direct monomial count gives
\begin{equation*}
C_{n,d,k} = \sum_{r=0}^{d-1} (r+1) \binom{d+k-r+n-2}{n-2},
\end{equation*}
with the convention that the sum is zero when $d=0$.

There is a particularly useful expression for the complementary dimension.

\begin{lemma}
\label{lem:H-dimension}
One has
\begin{equation*}
N_D-C_{n,d,k} = (d+1)R_k+B_k.
\end{equation*}
Equivalently, \( \dim[I_{\cH}^d]_{d+k} = (d+1)R_k+B_k. \)
\end{lemma}

\begin{proof}
Apply Proposition~\ref{prop:sections-identification} with $Z=\varnothing$. Then \( [I_{\cH}^d]_{d+k} \simeq H^0(\PP^1,\cT_k(d)). \) By Proposition~\ref{prop:Tk-splitting},
\[
\cT_k \simeq \bigoplus_{j=0}^k \cO_{\PP^1}(j)^{ \oplus\binom{n+k-j-2}{n-2}}.
\]
Since $d,j\geq0$, we have \( h^0\bigl(\PP^1,\cO_{\PP^1}(d+j)\bigr) = d+j+1. \) Therefore
\begin{equation*}
\dim[I_{\cH}^d]_{d+k} = \sum_{j=0}^k (d+j+1) \binom{n+k-j-2}{n-2} = (d+1)R_k+B_k.
\end{equation*}
This proves the assertion.
\end{proof}

For the fat point scheme $Z$, define
\begin{equation*}
\eta_Z(D) := h^1\bigl(\PP^n,\cI_Z(D)\bigr).
\end{equation*}
From
\[
0 \longrightarrow \cI_Z(D) \longrightarrow \cO_{\PP^n}(D) \longrightarrow \cO_Z \longrightarrow 0
\]
we obtain
\begin{equation}
\label{eq:h0IZ}
h^0\bigl(\PP^n,\cI_Z(D)\bigr) = N_D-\deg Z+\eta_Z(D).
\end{equation}

There are two related virtual dimensions which will be useful.

\begin{definition}
\label{def:virtual-dimensions}
The \emph{naive virtual dimension} is
\begin{equation*}
\vdim_{\mathrm{naive}} := N_D-\deg Z-C_{n,d,k}.
\end{equation*}

The \emph{base virtual dimension} is
\begin{equation*}
\vdim_{\mathrm{base}} := h^0\bigl(\PP^n,\cI_Z(D)\bigr) - C_{n,d,k}.
\end{equation*}
\end{definition}

Thus the naive virtual dimension treats the fat scheme $Z$ as imposing $\deg Z$ independent conditions in degree $D$, whereas the base virtual dimension takes into account the actual dimension of the degree-$D$ linear system through $Z$.

By \eqref{eq:h0IZ},
\begin{equation}
\label{eq:vdim-relation}
\vdim_{\mathrm{base}} = \vdim_{\mathrm{naive}} + \eta_Z(D).
\end{equation}

The expected dimension relevant to unexpectedness is therefore
\begin{equation*}
\edim := \max\{0,\vdim_{\mathrm{base}}\}.
\end{equation*}

Accordingly, the system \( [I_Z\cap I_{\cH}^d]_{d+k} \) is unexpected if \( \dim[I_Z\cap I_{\cH}^d]_{d+k} > \edim. \)

\subsection{The excess formula}

The failure of the splitting exponents to lie below $d+1$ will play a distinguished role.

\begin{definition}
\label{def:splitting-tail}
For fixed $d$ and $k$, define the \emph{splitting tail}
\begin{equation*}
T_{d,k}(Z,\cH) := \sum_{i=1}^{R_k} \max\{0,a_i-d-1\}.
\end{equation*}
\end{definition}

When no confusion is possible, we simply write $T$.

\begin{theorem}[Excess formula]
\label{thm:excess-formula}
Let \( Z=m_1P_1+\cdots+m_sP_s\subset\PP^n \) be a fat point scheme and let $\cH$ be a general codimension-two linear subspace. Put $D=d+k$. Then
\begin{equation}
\label{eq:adim-expanded}
\dim[I_Z\cap I_{\cH}^d]_{d+k} = (d+1)R_k+B_k-\Lambda_k(Z)+T.
\end{equation}

Moreover, \( \dim[I_Z\cap I_{\cH}^d]_{d+k} - \vdim_{\mathrm{naive}} = T+\delta_k(Z), \) while
\begin{equation}
\label{eq:excess-base}
\dim[I_Z\cap I_{\cH}^d]_{d+k} - \vdim_{\mathrm{base}} = T+\delta_k(Z)-\eta_Z(D).
\end{equation}
\end{theorem}

\begin{proof}
For every integer $a$ one has
\[
\max\{0,d-a+1\} = d+1-a+\max\{0,a-d-1\}.
\]
Applying this identity to Theorem~\ref{thm:dimension-splitting} gives
\begin{equation*}
\dim[I_Z\cap I_{\cH}^d]_{d+k} = \sum_{i=1}^{R_k} \left( d+1-a_i+\max\{0,a_i-d-1\} \right) = (d+1)R_k - \sum_{i=1}^{R_k}a_i + T.
\end{equation*}
Using Proposition~\ref{prop:sum-exponents}, \( \sum_{i=1}^{R_k}a_i = \Lambda_k(Z)-B_k, \) and therefore
\[
\dim[I_Z\cap I_{\cH}^d]_{d+k} = (d+1)R_k+B_k-\Lambda_k(Z)+T.
\]

On the other hand, Lemma~\ref{lem:H-dimension} gives
\[
\begin{aligned}
\vdim_{\mathrm{naive}} = N_D-\deg Z-C_{n,d,k} = (d+1)R_k+B_k-\deg Z.
\end{aligned}
\]
Since \( \deg Z = \Lambda_k(Z)+\delta_k(Z), \) we obtain
\[
\vdim_{\mathrm{naive}} = (d+1)R_k+B_k - \Lambda_k(Z) - \delta_k(Z).
\]
Subtracting this expression from \eqref{eq:adim-expanded} yields \( \dim[I_Z\cap I_{\cH}^d]_{d+k} - \vdim_{\mathrm{naive}} = T+\delta_k(Z). \)

Finally, by \eqref{eq:vdim-relation}, \( \vdim_{\mathrm{base}} = \vdim_{\mathrm{naive}}+\eta_Z(D), \) and hence
\[
\dim[I_Z\cap I_{\cH}^d]_{d+k} - \vdim_{\mathrm{base}} = T+\delta_k(Z)-\eta_Z(D).
\]
\end{proof}

The three terms in \eqref{eq:excess-base} have different geometric origins. The splitting tail $T$ measures the contribution coming from large splitting exponents. The term $\delta_k(Z)$ measures the part of the fat point structure which is invisible to the degree-$k$ fibrewise sections. Finally, $\eta_Z(D)$ measures the failure of $Z$ itself to impose independent conditions on hypersurfaces of degree $D$.

The distinction between the virtual and expected dimensions is important when $\vdim_{\mathrm{base}}<0$.

\begin{corollary}
\label{cor:unexpected-criterion}
If \( \vdim_{\mathrm{base}}\geq0, \) then \( [I_Z\cap I_{\cH}^d]_{d+k} \) is unexpected if and only if
\begin{equation*}
T+\delta_k(Z)>\eta_Z(D).
\end{equation*}
If \( \vdim_{\mathrm{base}}<0, \) then the expected dimension is zero and the system is unexpected if and only if
\[
\dim[I_Z\cap I_{\cH}^d]_{d+k}>0.
\]
\end{corollary}

\begin{remark}
\label{rem:raw-vs-expected-excess}
The quantity \( T+\delta_k(Z)-\eta_Z(D) \) is the difference between the actual dimension and the \emph{base virtual dimension}. It should not be confused with the difference between the actual and expected dimensions when $\vdim_{\mathrm{base}}<0$. In that case, \( \edim=0, \) so the unexpected excess is simply the actual dimension.
\end{remark}

There are two useful special cases.

\begin{corollary}
\label{cor:no-jet-defect}
Suppose that \( m_i\leq k+1\) for every \(i.\) Then \( \delta_k(Z)=0 \) and consequently
\begin{equation*}
\dim[I_Z\cap I_{\cH}^d]_{d+k} - \vdim_{\mathrm{base}} = T-\eta_Z(D).
\end{equation*}
\end{corollary}

\begin{corollary}
\label{cor:reduced-case}
If $Z$ is reduced, then \( \Lambda_k(Z)=\deg Z=|Z|,\; \delta_k(Z)=0, \) and
\begin{equation*}
\sum_{i=1}^{R_k}a_i = |Z|-B_k.
\end{equation*}
Moreover, \( \dim[I_Z\cap I_{\cH}^d]_{d+k} - \vdim_{\mathrm{base}} = T-\eta_Z(D). \)
\end{corollary}

\section{Reduction of excessive multiplicities}
\label{sec:reduction}

The effective-length computation suggests that multiplicities larger than $k$ should have a particularly simple effect on the bundle $\cE_{k,Z,\cH}$. We show that this is indeed the case: every multiplicity exceeding $k$ produces only a twist by a divisor on $\PP^1$. Geometrically, this twist corresponds to a fixed hyperplane component.

Let \( Z=m_1P_1+\cdots+m_sP_s\subset\PP^n \) and let $\cH$ be a general codimension-two linear subspace disjoint from $\Supp(Z)$. As before, we put \( q_i=\rho(P_i)\in\PP^1. \) We assume that the points $q_i$ are pairwise distinct, which holds for a general choice of $\cH$.

For a fixed $k\geq0$, set
\begin{equation}
\label{eq:ri-definition}
r_i:=(m_i-k)_+ = \max\{m_i-k,0\},
\end{equation}
and
\begin{equation*}
r:=\sum_{i=1}^s r_i.
\end{equation*}
We also introduce the truncated fat point scheme
\begin{equation}
\label{eq:Z-truncated}
Z_{\leq k} := \sum_{i=1}^s \min\{m_i,k\}P_i.
\end{equation}

Thus $Z_{\leq k}$ is obtained from $Z$ by replacing every multiplicity larger than $k$ by $k$.

\subsection{The local truncation phenomenon}

We begin with the local statement underlying the reduction.

\begin{lemma}
\label{lem:local-truncation}
Let $P\in\PP^n\setminus\cH$, let $q=\rho(P)$, and suppose that $m>k$. Then, in a neighbourhood of $q$,
\begin{equation*}
\cE_{k,mP,\cH} = \cE_{k,kP,\cH}\otimes\cO_{\PP^1}(-(m-k)q).
\end{equation*}
\end{lemma}

\begin{proof}
Choose a local parameter $t$ on $\PP^1$ centred at $q$ and local fibre coordinates \( u_1,\ldots,u_{n-1} \) centred at $P$. After completion, a local section of $\cT_k$ can be written uniquely in the form
\begin{equation}
\label{eq:local-section-Tk}
\sum_{|\beta|\leq k} c_\beta(t)u^\beta,
\end{equation}
where \( u^\beta = u_1^{\beta_1}\cdots u_{n-1}^{\beta_{n-1}} \) and \( c_\beta(t)\in\K[[t]]. \)

The condition that \eqref{eq:local-section-Tk} vanish to order at least $m$ at $P$ is equivalent to
\begin{equation}
\label{eq:ord-condition-m}
\operatorname{ord}_t c_\beta(t) \geq m-|\beta| \qquad \text{for every }|\beta|\leq k.
\end{equation}
Since $m>k$, the right-hand side is positive for every such $\beta$. We may therefore write \( c_\beta(t) = t^{m-k}c'_\beta(t). \) Condition \eqref{eq:ord-condition-m} becomes \( \operatorname{ord}_t c'_\beta(t) \geq k-|\beta|. \) But this is precisely the condition that \( \sum_{|\beta|\leq k} c'_\beta(t)u^\beta \) vanish to order at least $k$ at $P$.

Thus the local module defining $\cE_{k,mP,\cH}$ is obtained from the one defining $\cE_{k,kP,\cH}$ by multiplication by $t^{m-k}$. Equivalently,
\[
\cE_{k,mP,\cH} = \cE_{k,kP,\cH}\otimes \cO_{\PP^1}(-(m-k)q)
\]
locally at $q$.

Since completion is faithfully flat, the equality of the completed stalks gives the desired equality of the corresponding local subsheaves.
\end{proof}

\subsection{The truncation theorem}

The local statement globalizes immediately.

\begin{theorem}[Truncation theorem]
\label{thm:truncation}
Let \( Z=m_1P_1+\cdots+m_sP_s \) and define $Z_{\leq k}$ and $r_i$ as in \eqref{eq:Z-truncated} and \eqref{eq:ri-definition}. Put \( \Delta_{k,Z} := \sum_{i=1}^s r_iq_i. \) Then
\begin{equation*}
\cE_{k,Z,\cH} \simeq \cE_{k,Z_{\leq k},\cH} \otimes \cO_{\PP^1}(-\Delta_{k,Z}).
\end{equation*}
In particular, since \( \deg\Delta_{k,Z}=r, \) one has
\begin{equation}
\label{eq:bundle-truncation}
\cE_{k,Z,\cH} \simeq \cE_{k,Z_{\leq k},\cH}(-r).
\end{equation}
\end{theorem}

\begin{proof}
Away from the points \( q_1,\ldots,q_s, \) the schemes $Z$ and $Z_{\leq k}$ impose no conditions on the fibres of $\cT_k$, so the corresponding subsheaves coincide.

At a point $q_i$, the conditions coming from all points $P_j$ with $j\neq i$ are absent because the $q_j$ are pairwise distinct. If $m_i\leq k$, the schemes $Z$ and $Z_{\leq k}$ impose the same local condition. If $m_i>k$, Lemma~\ref{lem:local-truncation} gives \( \cE_{k,Z,\cH} = \cE_{k,Z_{\leq k},\cH} \otimes \cO_{\PP^1}(-r_iq_i) \) locally at $q_i$.

Combining these local identities gives
\[
\cE_{k,Z,\cH} \simeq \cE_{k,Z_{\leq k},\cH} \otimes \cO_{\PP^1} \left( -\sum_{i=1}^s r_iq_i \right).
\]
Since every divisor of degree $r$ on $\PP^1$ determines a line bundle isomorphic to $\cO_{\PP^1}(r)$, this yields
\eqref{eq:bundle-truncation}.
\end{proof}

The effect on the splitting type is particularly simple.

\begin{corollary}
\label{cor:splitting-shift}
Suppose that \( \cE_{k,Z_{\leq k},\cH} \simeq \bigoplus_{j=1}^{R_k} \cO_{\PP^1}(-b_j). \) Then
\begin{equation*}
\cE_{k,Z,\cH} \simeq \bigoplus_{j=1}^{R_k} \cO_{\PP^1}(-(b_j+r)).
\end{equation*}
Consequently, the splitting exponents satisfy
\begin{equation*}
a_j(Z) = a_j(Z_{\leq k})+r, \qquad j=1,\ldots,R_k.
\end{equation*}
\end{corollary}

Thus multiplicities larger than $k$ do not change the relative position of the splitting exponents: they translate the entire splitting type by the same integer $r$.

\subsection{Effective length under truncation}

The truncation theorem has an exact numerical counterpart.

\begin{proposition}
\label{prop:lambda-truncation}
For every $m\geq1$,
\begin{equation}
\label{eq:lambda-truncation}
\lambda_{n,k}(m) = \binom{n+\min\{m,k\}-1}{n} + (m-k)_+R_k.
\end{equation}
Consequently, \( \Lambda_k(Z) = \deg Z_{\leq k} + rR_k. \)
\end{proposition}

\begin{proof}
If $m\leq k$, then \( \lambda_{n,k}(m) = \length(mP) = \binom{n+m-1}{n}, \) and there is nothing to prove.

Suppose therefore that $m>k$. By definition,
\[
\lambda_{n,k}(m) = \sum_{b=0}^k (m-b) \binom{n+b-2}{n-2}.
\]
Writing \( m-b=(m-k)+(k-b), \) we obtain
\begin{equation*}
\lambda_{n,k}(m) = (m-k) \sum_{b=0}^k \binom{n+b-2}{n-2} + \sum_{b=0}^k (k-b) \binom{n+b-2}{n-2}.
\end{equation*}
The first sum equals \( R_k = \binom{n+k-1}{n-1}, \) while the second equals \[ \lambda_{n,k}(k) = \length(kP) = \binom{n+k-1}{n}. \] Hence 
\[ \lambda_{n,k}(m) = (m-k)R_k + \binom{n+k-1}{n}, \] which is \eqref{eq:lambda-truncation}.

Summing over all points of $Z$ gives
\[
\Lambda_k(Z) = \deg Z_{\leq k} + \left( \sum_{i=1}^s r_i \right)R_k,
\]
and therefore \( \Lambda_k(Z) = \deg Z_{\leq k}+rR_k. \)
\end{proof}

\begin{remark}
\label{rem:delta-vs-truncation}
There are two different thresholds in the behaviour of a fat point. The jet defect \( \delta_k(mP)= \length(mP)-\lambda_{n,k}(m) \) is zero whenever \( m\leq k+1. \) On the other hand, the truncation phenomenon of Theorem~\ref{thm:truncation} already occurs when \( m>k. \) Thus a point of multiplicity $k+1$ has no jet defect, but it already forces one fixed hyperplane component in the geometric interpretation below.
\end{remark}

\subsection{Fixed hyperplane components}

For each $P_i\in\Supp(Z)$, let \( L_i:=\langle\cH,P_i\rangle\subset\PP^n. \) Since $\cH$ has codimension two and $P_i\notin\cH$, the linear space $L_i$ is a hyperplane containing $\cH$.

Let \( \ell_i=0 \) be an equation of $L_i$. Since $\cH$ is general and the points $q_i=\rho(P_i)$ are pairwise distinct, one has
\[
P_j\notin L_i \qquad \text{for }i\neq j.
\]

The bundle-theoretic truncation has the following concrete interpretation.

\begin{theorem}[Fixed component reduction]
\label{thm:fixed-component-reduction}
With the notation above, suppose that \( d\geq r. \) Then multiplication by \( M_Z := \prod_{i=1}^s \ell_i^{r_i} \) induces an isomorphism
\begin{equation*}
\left[ I_{Z_{\leq k}} \cap I_{\cH}^{\,d-r} \right]_{d-r+k} \xrightarrow{\ \cdot M_Z\ } \left[ I_Z \cap I_{\cH}^{\,d} \right]_{d+k}.
\end{equation*}
In particular, \( \prod_{i=1}^s \ell_i^{r_i} \) is a fixed factor of every hypersurface in \( \left[ I_Z\cap I_{\cH}^d \right]_{d+k}. \)
\end{theorem}

\begin{proof}
We first show that every element of the target is divisible by $M_Z$.

Fix an index $i$ with $m_i>k$. The local computation in Lemma~\ref{lem:local-truncation} shows that every section of $\cE_{k,Z,\cH}$ contains the factor corresponding to the local parameter at $q_i$ to order at least \( r_i=m_i-k. \) Under the projection from $\cH$, this local parameter is represented by the hyperplane $L_i$. Hence every hypersurface in \( [I_Z\cap I_{\cH}^d]_{d+k} \) is divisible by \( \ell_i^{r_i}. \) Since the hyperplanes $L_i$ are distinct, every such hypersurface is divisible by \( M_Z=\prod_i\ell_i^{r_i}. \)

Let therefore \( F=M_ZG. \) The degree of $G$ is \( (d+k)-r = (d-r)+k. \) Each $\ell_i$ belongs to $I_{\cH}\setminus I_{\cH}^2$. Hence $M_Z$ has $I_{\cH}$-adic order exactly $r$. 
Indeed, its initial form in
\[
\operatorname{gr}_{I_{\cH}}S \simeq (S/I_{\cH})[u,v]
\]
is a nonzero product of linear forms. Since $\operatorname{gr}_{I_{\cH}}S$ is a domain, the $I_{\cH}$-adic order is additive on products. Therefore, from \( F=M_ZG\in I_{\cH}^d \) it follows that \( G\in I_{\cH}^{d-r}. \)

At $P_i$, the factor $\ell_i^{r_i}$ has multiplicity $r_i$, while $\ell_j(P_i)\neq0$ for $j\neq i$. Thus division by $M_Z$ lowers the multiplicity at $P_i$ precisely by $r_i$. Consequently,
\[
\mult_{P_i}(G) \geq m_i-r_i = \min\{m_i,k\}.
\]
Hence \( G \in \left[ I_{Z_{\leq k}} \cap I_{\cH}^{\,d-r} \right]_{d-r+k}. \)

Conversely, if $G$ belongs to the latter space, multiplication by $M_Z$ raises the multiplicity along $\cH$ by $r$ and raises the multiplicity at $P_i$ by $r_i$. Therefore \( M_ZG \in \left[ I_Z\cap I_{\cH}^d \right]_{d+k}. \) This proves the asserted isomorphism.
\end{proof}

\begin{remark}
\label{rem:bundle-fixed-component}
Theorem~\ref{thm:fixed-component-reduction} is the geometric counterpart of the bundle identity \( \cE_{k,Z,\cH} \simeq \cE_{k,Z_{\leq k},\cH}(-r). \) Indeed, twisting by $d$ gives \( \cE_{k,Z,\cH}(d) \simeq \cE_{k,Z_{\leq k},\cH}(d-r), \) and therefore
\[
H^0\bigl( \PP^1,\cE_{k,Z,\cH}(d) \bigr) \simeq H^0\bigl( \PP^1, \cE_{k,Z_{\leq k},\cH}(d-r) \bigr).
\]
When $d\geq r$, Theorem~\ref{thm:fixed-component-reduction} identifies this isomorphism explicitly with division by the fixed product $M_Z$.
\end{remark}

\subsection{Consequences for small values of $k$}

The truncation theorem clarifies where genuinely new phenomena caused by a nonreduced fixed scheme can first occur.

\begin{corollary}
\label{cor:k1-reduction}
Let $k=1$. Then \( Z_{\leq1} = \Supp(Z)_{\mathrm{red}} \) and \( r = \sum_{i=1}^s(m_i-1). \)
Consequently,
\begin{equation*}
\cE_{1,Z,\cH} \simeq \cE_{1,\Supp(Z)_{\mathrm{red}},\cH}(-r).
\end{equation*}
In particular, for $d\geq r$, every linear system of type $(d+1,d)$ with fixed fat point scheme $Z$ is obtained from the corresponding system for the reduced support by multiplication by the fixed hyperplanes $L_i^{m_i-1}$.
\end{corollary}

Thus, from the point of view of the bundle and its splitting type, fat point schemes introduce no essentially new phenomenon for $k=1$.

For $k=2$, after truncation every fixed point has multiplicity at most two. Hence the first genuinely nonreduced case not explained by fixed hyperplane components is a scheme containing double points. This is precisely the situation that will occur in the planar example of Section~\ref{sec:examples}.

\section{Examples}
\label{sec:examples}

We now illustrate the preceding results on explicit fat point configurations. The first example is planar and shows that the passage from reduced to nonreduced fixed schemes produces genuinely new unexpected curves even when all fixed multiplicities are bounded by $k$. In particular, the phenomenon is not explained by the fixed components of Section~\ref{sec:reduction}.

All computer algebra computations in this section were performed with \textsc{Singular} \cite{Singular}.

\subsection{An unexpected sextic with three double base points}
\label{subsec:planar-example}

Consider the following set of thirteen points in $\PP^2$:
\[
\begin{aligned}
P_1&=(1:0:0),&
P_2&=(-1:1:0),&
P_3&=(0:1:0),&
P_4&=(1:1:0),\\
P_5&=(-1:-1:1),& P_6&=(1:-1:1),& P_7&=(-1:1:1),& P_8&=(1:1:1),\\
P_9&=(0:-1:1),& P_{10}&=(-1:0:1),& P_{11}&=(0:0:1),& P_{12}&=(1:0:1),\\
&&P_{13}&=(0:1:1).
\end{aligned}
\]

We replace three of these points by double points and consider the fat point scheme
\begin{equation}
\label{eq:planar-Z}
Z = 2P_1+2P_3+2P_{11} + \sum_{i\notin\{1,3,11\}}P_i.
\end{equation}
Thus \( \deg Z = 3\cdot3+10 = 19. \)

We take \( k=2. \) Since every multiplicity of $Z$ is at most $k$, one has \( Z_{\leq2}=Z \) and the reduction theorem of Section~\ref{sec:reduction} does not remove any fixed hyperplane component. Moreover,
\[ \Lambda_2(Z)=\deg Z=19, \qquad \delta_2(Z)=0.
\]
In this case \( R_2 = \binom{3}{1} = 3,\) \( B_2 = \binom{3}{2} = 3. \) Consequently, Proposition~\ref{prop:sum-exponents} gives
\[
\sum_{i=1}^3a_i = 19-3 = 16.
\]
An exact rank computation over $\Q(a,b)$, with
$Q=(a:b:1)$ general, gives
\[
\dim[I_Z\cap I_Q^4]_6=1, \qquad \dim[I_Z\cap I_Q^5]_7=2.
\]
Together with $\sum_i a_i=16$, these dimensions determine the splitting type uniquely:
\[
\cE_{2,Z,Q} \simeq \cO_{\PP^1}(-4) \oplus \cO_{\PP^1}(-6)^{\oplus2}.
\]
Hence there is a unique, up to scalar, sextic containing $Z$ and having multiplicity at least four at a general point $Q$.

On the other hand, a Hilbert function computation gives \( \eta_Z(6)=0. \) Therefore \( h^0(\PP^2,\cI_Z(6)) = \binom{8}{2}-19 = 9. \) A point of multiplicity four in $\PP^2$ imposes \( \binom{5}{2}=10 \) conditions on sextics. Thus
\[
\vdim_{\mathrm{base}} = 9-10 = -1,
\]
and consequently \( \edim=0. \) We have proved the following.

\begin{proposition}
\label{prop:unexpected-sextic}
The fat point scheme $Z$ defined in \eqref{eq:planar-Z} admits an unexpected sextic with a general quadruple point. More precisely, for a general point $Q\in\PP^2$,
\[\dim[I_Z\cap I_Q^4]_6=1,
\]
whereas the expected dimension is zero.
\end{proposition}

This example also gives a particularly transparent illustration of the excess formula. For $d=4$ and splitting type $(4,6,6)$, the splitting tail is
\[
T = (4-5)_+ + (6-5)_+ + (6-5)_+ = 2.
\]
Since \( \delta_2(Z)=0\) and \( \eta_Z(6)=0, \) Theorem~\ref{thm:excess-formula} gives
\[
\dim[I_Z\cap I_Q^4]_6 - \vdim_{\mathrm{base}} = 2,
\]
which is exactly \( 1-(-1)=2. \) Notice that the excess over the expected dimension is only \( 1-0=1, \) in accordance with Remark~\ref{rem:raw-vs-expected-excess}.

The splitting type also shows that $d=4$ is the only value relevant to unexpectedness for $k=2$. Indeed, for $d<4$ the system has no sections, whereas for $d\geq5$ the splitting tail vanishes.

\subsubsection*{Absolute irreducibility}

The unexpected curve in Proposition~\ref{prop:unexpected-sextic} need not arise from a reducible curve with forced components. For instance, take \( Q=(-31:-9:15). \) The one-dimensional space \( [I_Z\cap I_Q^4]_6 \) is generated by an absolutely irreducible sextic, as verified by an exact absolute factorization computation in \textsc{Singular}.

Thus the unexpectedness in this example is not caused by a forced factorization of the sextic. It occurs for a fat point scheme whose fixed multiplicities satisfy \( m_i\leq k, \) and the resulting unexpected curve can be absolutely irreducible. This distinguishes the example from the reductions described in Section~\ref{sec:reduction}.

\subsection{The same distinguished triple with multiplicity three}
\label{subsec:planar-triple-points}

The three points \( P_1, P_3, P_{11} \) appearing in Proposition~\ref{prop:unexpected-sextic} exhibit a second unexpected phenomenon when their multiplicities are increased once again.

Let \( X=P_1+\cdots+P_{13} \) denote the reduced configuration considered in Subsection~\ref{subsec:planar-example}, and put
\begin{equation*}
Z_3 = 3P_1+3P_3+3P_{11} + \sum_{i\notin\{1,3,11\}}P_i.
\end{equation*}
Thus $Z_3$ consists of three triple points and ten simple points, and \( \deg Z_3 = 3\binom{4}{2}+10 = 28. \)

We now take \( k=3. \) Again, all fixed multiplicities are bounded by $k$, so the reduction theorem of Section~\ref{sec:reduction} does not remove any fixed component. Moreover, \( \Lambda_3(Z_3)=\deg Z_3=28, \delta_3(Z_3)=0. \) Since $R_3=4$ and $B_3=6$, Proposition~\ref{prop:sum-exponents} gives \( \sum_{i=1}^4 a_i=28-6=22. \)

An exact rank computation over $\Q(a,b)$, with $Q=(a:b:1)$ general, gives
\[
\dim[I_{Z_3}\cap I_Q^4]_7=0, \qquad \dim[I_{Z_3}\cap I_Q^5]_8=3.
\]
Together with $\sum_i a_i=22$, these dimensions determine the splitting type uniquely:
\begin{equation}
\label{eq:planar-Z3-splitting}
\cE_{3,Z_3,Q} \simeq \cO_{\PP^1}(-5)^{\oplus3} \oplus \cO_{\PP^1}(-7).
\end{equation}
Thus the splitting type is \( (5,5,5,7). \) The dimensions reconstructed from the splitting type begin as follows:
\[
\begin{array}{c|cccccccc}
d   &0&1&2&3&4&5&6&7\\
\hline
h^0\bigl(\PP^1,\cE_{3,Z_3,Q}(d)\bigr)
   &0&0&0&0&0&3&6&10.
\end{array}
\]
In particular, for \( d=5, D=d+k=8, \) we obtain
\begin{equation}
\label{eq:planar-octic-adim}
\dim[I_{Z_3}\cap I_Q^5]_8 = 3.
\end{equation}

On the other hand, a direct Hilbert function computation gives \( \eta_{Z_3}(8)=0. \) Hence \( h^0(\PP^2,\cI_{Z_3}(8)) = \binom{10}{2}-28 = 17. \) A point of multiplicity five in $\PP^2$ imposes \( \binom{6}{2}=15 \) conditions. Therefore
\[
\vdim_{\mathrm{base}} = 17-15 = 2,
\]
and hence \( \edim=2. \) Together with \eqref{eq:planar-octic-adim}, this gives \( \dim[I_{Z_3}\cap I_Q^5]_8 = 3 > 2 = \edim. \)

Thus the same distinguished triple of points which produces the unexpected sextic for $k=2$ also produces an unexpected system for $k=3$, now of type \( (8,5). \) The excess formula is particularly simple in this case. For $d=5$, the splitting tail of $(5,5,5,7)$ equals
\[
T = 3(5-6)_+ + (7-6)_+ = 1.
\]
Since \( \delta_3(Z_3)=0\) and \( \eta_{Z_3}(8)=0, \) Theorem~\ref{thm:excess-formula} gives \( \dim[I_{Z_3}\cap I_Q^5]_8 - \vdim_{\mathrm{base}} = 1, \) in agreement with the direct computation.

\subsubsection*{Relation with the unexpected sextic}

There is a natural relation between this three-dimensional system of octics and the unexpected sextic from Proposition~\ref{prop:unexpected-sextic}.

Write
\[
Z_2 = 2P_1+2P_3+2P_{11} + \sum_{i\notin\{1,3,11\}}P_i,\]
so that $Z_2$ is the fat point scheme considered previously. For a general point $Q$, let \( F_6(Q) \) be the unique, up to scalar, element of \( [I_{Z_2}\cap I_Q^4]_6. \)

Let \( V_Q = [I_{\{P_1,P_3,P_{11},Q\}}]_2 \) be the vector space of conics passing through \( P_1,P_3, P_{11},Q. \) Since these four points impose independent conditions on conics, \( \dim V_Q=2. \) For every $C\in V_Q$ we have \( F_6(Q)C \in [I_{Z_3}\cap I_Q^5]_8. \) Indeed, multiplication by $C$ raises the multiplicity at each of $P_1,P_3,P_{11}$ from two to at least three, and raises the multiplicity at $Q$ from four to at least five.

Consequently,
\begin{equation}
\label{eq:octic-inherited-subspace}
F_6(Q)V_Q \subset [I_{Z_3}\cap I_Q^5]_8
\end{equation}
is a two-dimensional subspace.

This observation gives a geometric interpretation of the dimension count. The expected dimension \( \edim=2 \) is already accounted for by the reducible octics
\[
F_6(Q)C, \qquad C\in V_Q.
\]
The actual system has dimension three, so it contains one additional independent direction which is not obtained in this way.

\subsubsection*{Irreducible members}

The unexpected octic system is not explained by reducible members arising from the construction above. For instance, consider the point
\[
Q=(-59:-21:35), \qquad I_Q=(5y+3z,7x-8y+7z).
\]
A direct computation gives \( \dim[I_{Z_3}\cap I_Q^5]_8=3. \) One of the degree-eight generators of this space is absolutely irreducible, as verified by an exact absolute factorization computation in \texttt{Singular}. Thus the linear system \( [I_{Z_3}\cap I_Q^5]_8 \) contains absolutely irreducible members. In particular, such a member cannot belong to the two-dimensional subspace \( F_6(Q)V_Q \) described in \eqref{eq:octic-inherited-subspace}.

\subsection{The Fermat configuration in $\PP^3$}
\label{subsec:fermat-P3}

We next consider a three-dimensional Fermat configuration. Let $\varepsilon$ be a primitive third root of unity and let \( X\subset\PP^3 \) be the set consisting of the four coordinate points
\[
(1:0:0:0),\quad (0:1:0:0),\quad (0:0:1:0),\quad (0:0:0:1)
\]
together with the $27$ points \( (1:\varepsilon^a:\varepsilon^b:\varepsilon^c), \) where \(a,b,c\in\{0,1,2\}. \) Thus \( |X|=31. \)

The reduced configuration $X$ has already appeared in the study of unexpected hypersurfaces of type $(d+k,d)$ in \cite{JMTG}. Here we modify it by replacing one of the coordinate points by a triple point. Put \( P=(0:1:0:0) \) and define
\begin{equation*}
Z = 3P+\sum_{Q\in X\setminus\{P\}}Q.
\end{equation*}
Thus $Z$ consists of one triple point and thirty simple points, and \( \deg Z = \binom{5}{3}+30 = 40. \)

We take \( k=3. \) Since all multiplicities of $Z$ are at most $k$, the truncation theorem of Section~\ref{sec:reduction} does not remove any fixed hyperplane component. Moreover, \( \Lambda_3(Z)=\deg Z=40\) and \( \delta_3(Z)=0. \) In the present case, \( R_3 = \binom{5}{2} = 10, B_3 = \binom{5}{3} = 10. \) Hence
\[
\sum_{i=1}^{10}a_i = \Lambda_3(Z)-B_3 = 30.
\]
Since $Z$ is defined over $\Q$, an exact rank computation over $\Q(\alpha,\beta,\gamma,\delta)$ for a general line
\[
x=\alpha z+\beta w,  \qquad   y=\gamma z+\delta w
\]
gives \( h^0\bigl(\PP^1,\cE_{3,Z,\cH}(d)\bigr) = 0,1,3,10 \) for $d=0,1,2,3$, respectively. Together with \( \sum_{i=1}^{10}a_i=30, \) these values determine the splitting type uniquely:
\begin{equation}
\label{eq:Fermat-fat-splitting}
\cE_{3,Z,\cH} \simeq \cO_{\PP^1}(-1) \oplus \cO_{\PP^1}(-2) \oplus \cO_{\PP^1}(-3)^{\oplus5} \oplus \cO_{\PP^1}(-4)^{\oplus3}.
\end{equation}
Thus the splitting type is \( (1,2,3,3,3,3,3,4,4,4). \) This single splitting type produces two unexpected linear systems.

\subsubsection*{The unexpected quartic}

First take \( d=1. \) Then \( D=d+k=4, \) and Theorem~\ref{thm:dimension-splitting} gives
\[
\begin{aligned}
\dim[I_Z\cap I_{\cH}]_4 = h^0\bigl(\PP^1,\cE_{3,Z,\cH}(1)\bigr) = h^0\bigl(\PP^1,\cO_{\PP^1}\bigr)=1.
\end{aligned}
\]

On the other hand, a direct Hilbert function computation gives \( h^0(\PP^3,\cI_Z(4))=5. \) Since \( h^0(\PP^3,\cO_{\PP^3}(4)) = \binom{7}{3} = 35, \) we have \( \eta_Z(4) = 5-35+40 = 10. \) A line in $\PP^3$ imposes \( h^0(\PP^1,\cO_{\PP^1}(4)) = 5 \) conditions on quartic hypersurfaces. Hence \( \vdim_{\mathrm{base}} = 5-5 = 0. \) Therefore \( \dim[I_Z\cap I_{\cH}]_4=1>\edim=0. \) Thus $Z$ admits an unexpected quartic containing a general line.

This phenomenon was already observed in \cite{JMTG}. The present approach recovers it directly from the splitting type \eqref{eq:Fermat-fat-splitting}.

The excess formula makes the numerical mechanism explicit. For $d=1$, the splitting tail equals
\[
\begin{aligned}
T = (1-2)_+ +(2-2)_+ +5(3-2)_+ +3(4-2)_+ = 5+6 =11.
\end{aligned}
\]
Since \( \delta_3(Z)=0\) and \( \eta_Z(4)=10, \) Theorem~\ref{thm:excess-formula} gives
\[
\dim[I_Z\cap I_{\cH}]_4 - \vdim_{\mathrm{base}} = 11-10 = 1.
\]

\subsubsection*{Unexpected quintics}

The same splitting type yields another unexpected system which is not visible at the level of quartics. Take \( d=2. \) Then \( D=d+k=5. \) Twisting \eqref{eq:Fermat-fat-splitting} by $2$ gives
\[
\cE_{3,Z,\cH}(2) \simeq \cO_{\PP^1}(1) \oplus \cO_{\PP^1} \oplus \cO_{\PP^1}(-1)^{\oplus5} \oplus \cO_{\PP^1}(-2)^{\oplus3}.
\]
Consequently,
\begin{equation}
\label{eq:Fermat-quintic-adim}
\dim[I_Z\cap I_{\cH}^2]_5 = 2+1 = 3.
\end{equation}
In degree five, the Hilbert function of $Z$ gives \( \eta_Z(5)=1. \) Since \( h^0(\PP^3,\cO_{\PP^3}(5)) = \binom{8}{3} = 56, \) we obtain
\[
h^0(\PP^3,\cI_Z(5)) = 56-40+1 = 17.
\]
Vanishing to order at least two along a line imposes \( C_{3,2,3}=16 \) conditions on quintics. Hence \( \vdim_{\mathrm{base}} = 17-16 = 1. \) It follows from \eqref{eq:Fermat-quintic-adim} that
\[
\dim[I_Z\cap I_{\cH}^2]_5=3>1=\edim.
\]
Thus $Z$ admits unexpected quintic surfaces vanishing doubly along a general line.

Again, the excess formula explains the discrepancy. For $d=2$, only the three exponents equal to $4$ contribute to the splitting tail, and therefore \( T=3. \) Since \( \delta_3(Z)=0\) and \( \eta_Z(5)=1, \) we obtain
\[
\dim[I_Z\cap I_{\cH}^2]_5 - \vdim_{\mathrm{base}} = 3-1 = 2.
\]

The relation between the quartic and quintic systems is also worth noting. Let $F_4$ denote the unique quartic in \( [I_Z\cap I_{\cH}]_4. \) If \( I_{\cH}=(L_1,L_2), \) then \( L_1F_4,L_2F_4 \) belong to \( [I_Z\cap I_{\cH}^2]_5. \)
They span a two-dimensional subspace. Since the full quintic system has dimension three, there is one additional independent direction which is not obtained by multiplying the unexpected quartic by a linear form vanishing on $\cH$.

For instance, consider the line
\[
\cH = V(x-3y+4z-6w,\;x+y-7z+11w).
\]
For this line the degree-five space \( [I_Z\cap I_{\cH}^2]_5 \) is three-dimensional. If $F_1,F_2,F_3$ denote the three degree-five generators obtained by the computation, then $F_1+F_2+F_3$ is absolutely irreducible, as verified by an exact absolute factorization computation in \textsc{Singular}. Hence the unexpected quintic system contains absolutely irreducible members.

A short \textsc{Singular} script reproducing the explicit computations used in the irreducibility checks above is available at \cite{Code}.

\end{document}